\documentclass{article}
\newtheorem{Th}{Theorem}
\newtheorem{Cor}{Corollary}
\newtheorem{Prop}{Proposition}
\newtheorem{Lm}{Lemma}
\newtheorem{Def}{Definition}

\newtheorem{rem}{Remark}
\newcommand{\X}{{\cal X}}

\newcommand{\R}{{\bf R}}
\newcommand{\N}{{\bf N}}
\newcommand{\C}{{\bf C}}
\newcommand{\ns}{\mbox{ns}}

\renewcommand{\H}{{\cal H}}
\newcommand{\U}{{\cal U}}
\newcommand{\V}{{\cal V}}
\newcommand{\K}{{\cal K}}

\renewcommand{\l}{\lambda}
\newcommand{\f}{\varphi}
\newcommand{\Mu}{{\rm M}}
\renewcommand{\a}{\alpha}
\renewcommand{\b}{\beta}
\newcommand{\supp}{\mbox{supp}}
\renewcommand{\*}{\,^*\!}\renewcommand{\o}{\,^\circ\,\!\!}
\newtheorem{Remark}{Remark}
\author{L.Yu. Glebsky, \and E.I.Gordon\thanks{The work of the second author was supported in part by NSF Grant 97 0009}}
\title{On approximation of topological groups by finite algebraic systems}
\date{}
\usepackage{latexsym}
\begin{document}
\maketitle

\begin{abstract}
It is known that locally compact groups approximable by finite ones are 
unimodular, but this condition is not sufficient, for example, the simple Lie groups are not approximable by finite ones as topological groups.  
In this paper the approximations of locally compact groups by more general finite algebraic systems are investigated. It is  proved that the approximation of locally compact groups by finite semigroups is equivalent to approximation by finite groups and thus not all locally compact groups are approximable by finite semigroups. We prove that any locally compact group is approximable by finite left (right) quasigroups but the approximabilty of a locally compact group by finite quasigroups (latin squares) implies its unimodularity. 
The question if the unimodularity of a locally compact group implies its approximability by finite 
quasigroups is open. We prove only that the discrete groups 
are approximable by finite quasigroups.
\end{abstract}

\section{Introduction}

The notion of approximation of a topological group by finite ones was introduced by the second author (cf. the monograph \cite{Gor} and the bibliography there). It was investigated in details for the case of locally compact abelian (LCA) groups in \cite{Gor} and for the case of discrete groups in
\cite{VG}. The approximations of LCA groups were used in \cite{AGKh} for a construction of finite
dimensional approxiamtions of pseudodifferential operators. The approximations of discrete groups have some interesting applications in the ergodic theory of group actions \cite{VG},\cite{AGG} and in symbolic dynamics \cite{Gr}. The approximability of any LCA group by finite abelian groups is proved in \cite{Gor}, the approximability of a huge class of nilpotent Lie groups by finite nilpotent groups is proved in \cite{VG}. The class of discrete approximable groups is a proper extension of the class of locally residually finite groups; there exist some non-approximable groups:
the Baumslag --- Solitar groups, finitely presented infinite simple groups and some other \cite{VG}.
It was proved in \cite{Resv} that all approximable locally compact groups are unimodular (the left and right Haar measures coincide). This condition is not sufficient --- we have mentioned already that there exist non-approximable discrete
groups. It was proved in \cite{AGG} that the simple Lie groups are not approximable by finite groups
as topological groups (since these groups are locally residually finite they are approximable as 
discrete groups).     

In this paper we investigate the approximation of locally compact groups by more general universal
algebras with one binary operation. We prove that any locally compact group is approximable by finite left (right) quasigroups --- the algebras, that have left (right) division\footnote{See Definition 2(1)}. Then we 
consider approximations of more general topological algebras and
prove that if any locally compact left (right) quasigroup $A$ that has the operation of taking left (right) inverse element (satisfies the left (right) cancellation law\footnote{See Definition 2(4)}) is approximable by finite left (right) quasigroups then there exists a positive non-trivial linear functional $I$ on $\C_0(A)$ (the space of all continuous functions on $A$ with compact support) such that $I(f)\geq I(l_h(f))$ ($I(f)\geq I(r_h(f))$) for any 
non-negative $f\in C_0(A)$ and any
$h\in A$. 
Here $l_h$ ($r_h)$ is
the left (right) shift on $\C_0(A)$, i.e. $l_h(f)(a)=f(h\circ a)$ ($r_h(f)(a)=f(a\circ h)$, where $\circ$ is the operation in $A$.

These inequalities imply immediately that if $A$ is a group, then $I$ is left (right) invariant and thus if the group $A$ is approximable by finite quasigroups, i.e. algebraic systems that are left and right quasigroups simultaneously, then $A$ is unimodular. So, any locally compact group that 
is approximable by finite quasigroups is unimodular. It is an interesting open question if the 
approximability by finite quasigroups implies the unimodularity of a locally compact group. We may say only that
the class of locally compact groups that are approximable by finite quasigroups is larger than the class of locally compact
groups approximable by finite groups. We prove that all discrete groups are approximable by finite quasigroups.  

We prove also that the approximabilty of a locally compact group by finite semigroups implies its approximability by finite groups.

In the proofs of mentioned results we use the language of nonstandard analysis that allows to simplify essentially these proofs. The necessary notions and results of nonstandard analysis 
can be found in the monographs \cite{Alb} or \cite{Loeb}

The results of this paper were discussed at the seminar "Nonstandard analysis" in the UIUC. The authors are grateful to Prof. C.W.Henson and Prof. P.Loeb for many important remarks. 

\section{Formulation of the main results}
We mean by algebra here a universal algebra that contains only one binary operation.

Let $<A,\circ>$ be a locally compact Hausdorff topological universal algebra, $C\subseteq A$ --- a compact, $\U$ --- a finite cover of $C$ by
open sets, $<H,\odot>$ --- a finite universal algebra.
In what follows we'll omit the symbols of operation and denote these algebras by $A$ and 
$H$ respectively.

\begin{Def} \label{approx}
\begin{enumerate}
\item A set $M\subseteq A$ is called a ($C,\U$)-grid iff 
$$
\forall U\in{\U}\left((C\cap U\neq\emptyset)\Rightarrow
\exists m\in M\;(m\in U)\right).
$$
\item A map $j:H\to A$ is called a ($C, \U$)-homomorphism
iff
$$
\forall x,y\in H\left((j(x),j(y),j(x)\circ j(y)\in C)\Rightarrow 
\exists U\in{\U} \;\;(j(x\odot y)\in U \land j(x)\circ j(y)\in
U)\right). 
$$
\item We say that the pair $<H,j>$ is a $(C,\U)$-approximation of $G$ if $j(H)$ is a $(C,\U)$-grid and
$j:H\to G$ is a $(C,\U)$-homomorphism.
\item Let $\K$ be a class of finite algebras. We say that $A$ is approximable by the systems of the class $\K$ if for any compact $C\subseteq A$ and for any finite cover $\U$ of $C$ by open sets there exists a ($C,\U$)-approximation $<H,j>$ of $A$ such that $H\in\K$ and $j$ is an injection.
\end{enumerate}
\end{Def}

\begin{Remark}
It is easy to see that the similar definition can be formulated for any topological universal algebra and it is not necessary to assume that 
approximated 
algebras are finite. For example, the approximations of discrete groups by amenable ones were introduced in \cite{AGG}. The approximations 
of universal algebras with finite signatures will be considered in another paper.
\end{Remark}

It is easy to see that the following propositions hold.

\begin{Prop} \label{sep_approx}
If $A$ is separable as a topological space then $A$ is approximable by algebras of a class $\K$
iff there exist a sequence of finite algebras $<H_n,\circ_n>$, $H_n\in\K$ and a sequence of injections 
$j_n:H_n\to A$ such that
for any compact $C\subseteq A$ and for any finite cover ${\cal U}$ of $C$ by open sets
there exists an $n_0\in\N$, such that for any $n>n_0$ $<H_n,j_n>$ 
 is a ($C,\cal U$)-approximation of $A$. 
\end{Prop}

\begin{Prop} \label{dis_approx}
A discrete algebra $A$ is approximable by algebras of a class $\K$ iff for any finite subset
$S\subseteq A$ there exist an algebra $H\in\K$ and an injection $j:S\to H$ such that 
$$
\forall s_1,s_2\in S\left(s_1\circ s_2\in S\Rightarrow j(s_1\circ s_2)=j(s_1)\odot j_(s_2)\right)  
$$
\end{Prop}
 
The definition of approximation of a LC group by finite algebras can 
be simplified a little. 
Let $G$ be a LC group.
We will denote by $\cdot$ the multiplication in $G$ and use the usual notations
$$
XY=\{x\cdot y\;|\;x\in X,\,y\in Y\}
$$
$$
X^{-1}=\{x^{-1}\;|\;x\in X\}
$$
$$
gX=\{g\cdot x\;|\;x\in X\}
$$
for $X,Y\subseteq G,\ g\in G$.

\begin{Def} \label{gr_approx}
Let $C\subseteq G$ be a compact, $U$ --- a relatively compact neighborhood of 
the
unit in $G$, and $H$ --- a finite algebra.

\begin{enumerate}
\item We say that a set $M\subseteq G$ is 
an $U$-grid of $C$ 
iff $C\subseteq MU$.
\item A map $j:H\to G$ is called a $(C,U)$-homomorphism if
$$
\forall x,y\in H\;\left( (j(x),j(y),j(x)\cdot j(y)\in C)\Rightarrow
(j(x\odot y)\in j(x)j(y)U) \right)
$$
\item We say that the pair $<H,j>$ is a $(C,U)$-approximation of $G$ if 
$j(H)$ is 
an $U$-grid of $C$ 
and
$j:H\to G$ is a $(C,U)$-homomorphism.
\item Let $\K$ be a class of finite algebras. We say that $G$ is approximable by the systems of the class $\K$ if for any compact $C\subseteq G$ and for any neighborhood of 
the
unit $U$ there exists a ($C,U$)-approximation $<H,j>$ of $G$ such that $H\in\K$ and $j$ is an injection.
\end{enumerate}
\end{Def}

\begin{Prop} \label{equiv-for-groups}
A locally compact group $G$ is approximable by the systems of a class $\K$ in the sense of
Definition~\ref{gr_approx} iff it is approximable by the systems of $\K$ in the sense of Definition
\ref{approx}
\end{Prop}
This proposition will be proved in Section 4.

\begin{Def} \label{quasigroups}
\begin{enumerate}
\item We say that algebra $A$ is an l-quasigroup 
(r-quasigroup) iff for all $a,b\in A$ the
equation $a\circ x=b$ ($x\circ a=b$) has the unique solution
$x=/(b,a)$ ($x=\backslash (b,a)$).
\item If the functions $\circ (\cdot,\cdot):A^2\to A$
$/(\cdot,\cdot):A^2\to A$, ($\backslash (\cdot,\cdot):A^2\to A$) are continuous
then we say that $A$ is a topological l-quasigroup (r-quasigroup).
\item We say that $(A,\circ)$ is a (topological) quasigroup iff
it is a (topological) l- and r-quasigroup simultaneously.
\item We say that l-quasigroup (r-quasigroup) $A$ satisfies the l-cancellation law
(r-cancellation  law) iff  there exists a function $b:A\to A$
such that
$$
\forall a \forall x\, b(a)\circ(a\circ x)=x\;\;\;
(\forall a \forall x\, (x\circ a)\circ b(a) =x).
$$
We write $b_a$ instead of $b(a)$.
\end{enumerate}
\end{Def}

There is a huge literature concerning quasigroups, cf., for example, \cite{quas}.
The operation table of a finite quasigroup is a latin square --- an $n\times n$-table
of $n$ elements $\{a_1,\dots,a_n\}$ such that all elements in each row and in each column are
distinct. An $n\times n$-table with this property that contains more than $n$ elements is called a latin subsquare. It is known \cite{Ryser} that any $n\times n$ latin subsquare 
with $k$ distinct elements
can be completed to
an 
$r\times r$ latin square, where $r=\max\{2n,k\}$.
 This fact together with Proposition~\ref{dis_approx} imply
immediately the following 
\begin{Prop} \label{dis_quas_appr}
Any discrete quasigroup, and thus any discrete group, is approximable by finite quasigroups.
\end{Prop}  

\begin{Th} \label{groupappr}
Any locally compact group $G$ is approximable by finite $l$-quasigroups ($r$-quasigroups).
\end{Th}

This theorem will be proved in Section 3.

\begin{Th} \label{functional}
Let $A$ be a locally compact $l$-quasigroup ($r$-quasigroup) that satisfies the $l$-cancellation ($r$-cancellation) law and that is approximable by finite $l$-quasigroups ($r$-quasigroups). Then there
exists a positive bounded non-trivial linear functional $I$ on $\C_0(A)$ that satisfies the inequality
$$
I(f)\geq I(l_h(f))\quad \left(I(f)\geq I(r_h(f))\right) \eqno (1)
$$
for any non-negative $f\in C_0(A)$ and any $h\in A$.

If $A$ is a locally compact quasigroup that satisfies 
both cancellation laws and is approximable by finite quasigroups, 
then $I$ satisfies both inequalities (1) simultaneously.
\end{Th}

This theorem will be proved in Section 5. 

\begin{Prop} \label{invariant}
If $A$ is a group and $I$ is a functional on $\C_0(A)$ that satisfies (1) then $I$ is left (right)
invariant.
\end{Prop}

{\bf Proof}
$$
I(f)\geq I(l_h(f))\geq I\left(l_{h^{-1}}\left(l_h(f)\right)\right)=I(f)\quad\Box
$$

\begin{Remark} Theorem~\ref{groupappr}, \ref{functional} and Proposition~\ref{invariant}
imply the existence of the Haar measure on a locally compact group. 
Indeed, the proof of Theorem~\ref{groupappr} given in Section~3 includes 
some ideas and constructions of the proof of existence of the Haar measure 
contained in the famous monograph \cite{Halmos}.
Yet another proof of existence of the Haar measure based 
on nonstandard analysis is given in 
\cite{Ross}

\end{Remark}

\begin{Cor} Any locally compact group $G$ approximable by finite quasigroups is unimodular.
\end{Cor}

This corollary generalizes the result of \cite{Resv} that any locally compact group approximable by finite groups is unimodular. Proposition~\ref{dis_quas_appr} shows that the class of groups 
approximable by finite quasigroups is larger then the class of groups approximable by finite groups
since there exist discrete groups that are not approximabale by finite groups \cite{VG}

{\bf Conjecture} A locally compact group is unimodular iff it is approximable by finite quasigroups.

\begin{Th} \label{semigroups}
A locally compact group is approximable by finite semigroups iff it is approximable by finite groups.\end{Th}

This theorem will be proved in Section 6.

\section{Proof of Theorem~\ref{groupappr}}

In this section $G$ is a locally compact group.

To prove Theorem~\ref{groupappr} we introduce the following construction. 
Given a 
neighborhood of the unity $U$ and a compact $C$ we find a finite $U$- grid 
of $C$
(Definition~\ref{gr_approx} (1)) $F\subseteq G$ and a collection $\{A_{g,h}\subseteq F\;;\;g,h\in F\}$, 
such that the following lemma holds.

\begin{Lm} \label{Lm_l-quasi}
For any
neighborhood of the unity $U$ and any compact $C$ there exist a finite
$F\subseteq G$, and a collection 
$\{A_{g,h}\subseteq F\;;\;g,h\in F\}$,
satisfying the following
conditions:
\begin{enumerate}
\item $F$ is an $U$-grid of $C$;
\item if $g,h\in C\cap F$ then $A_{g,h}\subseteq ghU$;
\item 
$\forall g\in F \forall S\subseteq F\;\;|\bigcup_{h\in S}A_{g,h}|
\geq |S|$. 
\end{enumerate}
\end{Lm}

We also need the well known combinatorial Theorem of P. Hall (the Marriage Lemma), see, for example, 
\cite{Ryser}.
\begin{Def}\label{Def_SDR}
Let $F_i\subseteq F,\ i=1,\dots,m$.  We say that the sequence $F_1,F_2,...,F_m$ has a
system of distinct representatives (SDR) iff we can chose $m$-permutation
of $F$ $a_1, a_2, ..., a_m$, such that $a_i\in F_i$ for $i=1,...,m$.
(Definition of $m$-permutation requires that $a_i\neq a_j$ for $i\neq j$.)
\end{Def}

\begin{Th} \label{Th_Hall}
The subsets $F_1, F_2, ... F_m$ have an SDR if and only if for each 
$S\subseteq \{1,2,...,m\}$ the following inequality holds
$$
|\bigcup_{k\in S} F_k|\geq |S|
$$
\end{Th}

\begin{rem} Nonstandard analysis versions of P.Hall's theorem were investigated in \cite{Bos}
\end{rem}
Lemma~\ref{Lm_l-quasi} (3) and Theorem~\ref{Th_Hall} imply that the set $F$ may be
equipped with operation $\odot$, satisfying definition of l-quasigroup. 
Indeed, by the condition 3)
the system $\{A_{g,h}\ |\ h\in F\}$ satisfies Theorem~\ref{Th_Hall} for any 
fixed 
$g\in F$ 
and thus
for any 
$g,h\in F$ 
there exists 
$g\odot h\in A_{g,h}$ 
such that for any 
$g\in F$ 
$\{g\odot h\ |\ h\in F\}$ is a permutation of $F$. 
Thus $<F,\odot>$ is an $l$-quasigroup. The conditions 1) and 2) of Lemma~\ref{Lm_l-quasi} imply that the $l$-quasigroup $<F,\odot>$ 
with identical  inclusion
is a $(C,U)$-approximation of $G$, see Definition~\ref{gr_approx}(3).

So, to complete the proof of Theorem~\ref{groupappr} we have only to prove  
Lemma~\ref{Lm_l-quasi}.

Let $O\subseteq G$ be a neighborhood of the unity and $A\subseteq G$. 
Denote $(A:O)$ the minimal number $n$, such that there exist
$F$, $|F|=n$, $A\subseteq FO$. 

In the following Propositions~\ref{Prop_1.1} and \ref{Prop_1.2} we assume that
\begin{itemize}
\item $O$ is a neighborhood of the unity;
\item $K$ is a compact;
\item $F$ is a finite set, $|F|=(K:0)$ and $K\subseteq FO$
($F$ is an optimal $O$-grid of $K$).
\end{itemize}
\begin{Prop}\label{Prop_1.1}
Let $S\subseteq F$, then $(SO:O)=|S|$.
\end{Prop}
{\bf Proof.} It is clear that $(SO:O)\leq |S|$. Suppose, that
$(SO:O)< |S|$, then 
$$
K\subseteq SO\cup\bigcup_{x\in F\setminus S}xO,
$$
and we can cover $K$ with less then $|F|$ elements.  $\Box$

\begin{Prop}\label{Prop_1.2}  
Let $M\subseteq K$. Then $|MO^{-1}\cap F|\geq (M:O)$.
\end{Prop}
{\bf Proof}
One has $M\subseteq K\subseteq FO$. It means that 
$\forall x\in M\exists f\in F\exists \epsilon\in O\;\;x=f\epsilon$,
or $f=x\epsilon^{-1}$ 
so, $f\in MO^{-1}$. 
Consequently, $MO^{-1}\cap F$ is an $O$-grid of
$M$. So, 
$(M:O)\leq |MO^{-1}\cap F|$. 
$\Box$  

{\bf Proof of Lemma~\ref{Lm_l-quasi}.}
Given a neighborhood of the unity $U\subseteq G$ and a compact $C\subseteq G$
one can chose a neighborhood of the unity $O$ and a compact $K$, such that 
\begin{itemize}
\item $OO^{-1}\subseteq U$;
\item $C^2\subseteq K$;
\item $CU\subseteq K$.
\end{itemize}
Let $F$ be an optimal $O$-grid of $K$. Define the sets $A_{g,h}$ as follows:
$$
A_{g,h}=\left\{\begin{array}{ll}
                ghOO^{-1}\cap F,&\mbox{if}\;\;g,h\in C\\
                F,& \mbox{otherwise}.
                \end{array} 
         \right.       
$$
It is easy to see that $F$ is $U$-grid of $C$  and item~2) of the 
lemma~\ref{Lm_l-quasi} is also satisfied. Prove item~3). 
Nontrivial case is when
$g\in C$ and $S\subseteq C$. 
By Proposition~\ref{Prop_1.1} $(SO:O)=|S|$, consequently, $(gSO:O)=|S|$.
Then, by Proposition~\ref{Prop_1.2}, 
$$
|S|\leq |gSOO^{-1}\cap F|=|\bigcup_{h\in S}A_{g,h}|.\ \ \Box
$$ 

\section{Nonstandard analysis approach to approximation of algebras}

In this section we introduce a brief exposition of nonstandard analysis (see, for example, 
\cite{Gor},\cite{Alb} or \cite{Loeb} for details).

Let $\l$ be an infinite cardinal.
Consider a $\l^+$-saturated nonstandard extension $\*\V$ of the standard universe $\V$. Recall that $\*\V$ is $\l^+$-saturated if for any
family $\X$ of {\it internal} sets (i.e. the sets that are elements of 
$\*\V$), 
such that $|\X|\leq\l$
and $\X$ has the finite intersection property, 
holds $\bigcap\limits_{X\in\X}X\neq\emptyset$.

There exits an embedding $\*:\V\to\*\V$ that satisfies the {\it Transfer principle}. The image of
an element $v\in\V$ under this embedding is denoted by $\*v$ and is called the nonstandard extension of $v$. A proposition $\f(v_1,\dots,v_n)$  
about $v_1,\dots,v_n$
is called {\it internal} if it is a statement formulated in usual 
mathematical terms --- not including such notions as 
"standard element", "nonstandard extension", etc. More formally $\f$ can be a formula of the language of 
the
set theory, or of the language of the theory of superstructures, 
or of the language of the elementary analysis \cite{Gor}, etc., depending on what kind of standard universe $\V$ we consider.

{\bf Transfer principle} If $\f$ is an internal formula and $v_1,\dots,v_n\in\V$ then 
$\f(v_1,\dots,v_n)$ holds in $\V$ iff $\f(\*v_1,\dots,\*v_n)$ holds in $\*\V$.

{\bf Example} Let $n\in\N$ and $B_n=\{\xi\in\R\ |\ \xi>n\}$. By transfer principle $\*B_n=
\{\xi\in\*\R\ |\ \xi>\*n\}$ and $\*n\in\*\N$. The elements of $\*\N$ are called the hypernatural numbers and the elements of $\*\R$ --- the hyperreals. Usually the notation $\*B$ is used only for 
the nonstandard extension of a standard set $B$. The nonstandard extension of a standard element
$b$ is denoted by $b$ also. So, if $\xi\in B_n$ then by the Transfer principle $\xi\in \*B_n$ and thus $B_n\subseteq\*B_n$.

The countable sequence $\{\*B_n\ |\ n\in\N\}$ has obviously the finite intersection property
and thus by saturation $\Mu=\bigcap\limits_n\*B_n\neq\emptyset$. 
Let $-\Mu=\{-x\;:\;x\in\Mu\}$.
The elements of $\Mu\cup -\Mu$ are called the {\it infinite} elements of $\*\R$ since if $\eta\in\Mu$ then 
$$
\forall\xi\in\R (|\eta|>|\xi|) \eqno(2)
$$
The elements of $\*\R$ inverse to infinite elements (and $\*\R$ is an ordered field by Transfer principle) are called the {\it infinitesimals}. It follows from (2) that if $\a$ is an infinitesimal
then 
$$
\forall\xi\in\R(\xi>0\Rightarrow |\a|<\xi).
$$

The set of all infinitesimals is called {\it the monad of zero} and denoted by $\mu(0)$. The sets 
$\Mu$ and $\mu(0)$ are not internal --- such sets are called {\it external}. Indeed if $\mu(0)$ would be
internal, then, being bounded from above, it must have the supremum by the Transfer principle. But
it is easy to see that both conjectures: 1) $\sup\mu(0)$ is an infinitesimal and 2) $\sup\mu(0)$ is not infinitesimal --- lead to the contradiction.

Two elements $\xi, \eta\in\*\R$ are {\it infinitely close} ($\xi\approx\eta$) if $\xi-\eta\in\mu(0)$.
In particular, $\xi\in\mu(0)$ iff $\xi\approx 0$.

The elements of $\*\R$ that are not infinite are called {\it bounded} or {\it finite}. It can be proved that any bounded element $\xi$ is infinitely close to the unique standard element $\a$.
This $\a$ is called {\it the standard part} of $\xi$ or {\it the shadow} of $\xi$ and denoted by
$\o\xi$. Now it is easy to prove that any hypernatural number that is not infinite is standard.

An internal set $B$ is called {\it hyperfinite} if there exists an $n\in\*\N\setminus\N$ and an
internal bijection $\f:\{0,1,\dots,n-1\}\to B$. Then we say that $n$ is the cardinality of $B$ and
write $|B|=n$. When we deal with internal sets we use the term "cardinality" for internal cardinality.  By the Transfer principle the hyperfinite sets have many features of standard finite sets. For example, an internal subset of a hyperfinite set is hyperfinite itself. Any hyperfinite set $B$ of reals has the maximal and the minimal elements and the $\sum\limits_{b\in B}b$ is defined. If $B=B_1\cup B_2$, where $B_1$ and $B_2$ are internal sets and $B_1\cap B_2=\emptyset$ then
$$
\sum\limits_{b\in B}b=\sum\limits_{b\in B_1}b+\sum\limits_{b\in B_2}b.
$$

We consider now the nonstandard extension $\*X$ of a standard 
Hausdorff topological space $X$ of the weight
(the cardinality of a minimal base of topology on $X$) less or equal to $\l$. If $x\in X$ is not an isolated point of $X$ then for any open $U\ni x$ the set $U\setminus\{x\}\neq\emptyset$. Let
$$
\mu(x)=\bigcap\{\*U\ |U\ni x,\,U\ \mbox{is open}\}.
$$
The set $\mu(x)$ (external, if $x$ is not an isolated point) is called {\it the monad} of $x$. The 
$\l^+$-saturation of the nonstandard universe $\*\V$ implies immediately that $\mu(x)$ contains some
nonstandard elements of $\*X$ if $x$ is not an isolated point. We say that an element $y\in\*X$ is nearstandard if $y\in\mu(x)$ for some standard $x\in X$. The set of all nearstandard elements of
$\*X$ is denoted by $\ns(\*X)$. 
Obviously $\ns(\*X)=\bigcup\limits_{x\in X}\mu(x)$. 
If $y\in\mu(x)$ then we say as before that $x$ is the standard part or 
the shadow of $y$ and write $x=\o y$ (it is unique for a Hausdorff space). 
We say that two nearstandard elements $y_1$ and $y_2$ are infinitely close ($y_1\approx y_2$) if $\o y_1=\o y_2$.

We list the well known properties (\cite{Gor},\cite{Alb},\cite{Loeb}) of $\*X$ in the following
\begin{Prop}   \label{topprop}
\begin{enumerate}
\item A set $U\subseteq X$ is open iff $\forall x\in U\left(\mu(x)\subseteq\*U\right)$.
\item A set $F\subseteq X$ is closed iff $\forall x\in X\left(\mu(x)\cap\*F\neq\emptyset\Rightarrow
x\in F\right)$.
\item $X$ is compact iff $\ns(\*X)=\*X$.
\item $X$ is locally compact iff $\ns(\*X)=\bigcup\{\*C\ |\ C\subseteq X,\ C\ \mbox{is compact}\}$.
\item If $Y$ is a topological space of the weight less or equal to $\l$ 
then a mapping $f:X\to Y$
is continuous iff $\forall x\in X,\xi\in\*X\left(\xi\approx x\Longrightarrow \*f(\xi)\approx f(x)\right)$. 
\end{enumerate}
\end{Prop} 

Let us return now to a locally compact algebra $A$. Assume that the cardinality of its topology (the family of all open sets) is $\l$ and consider the nonstandard extension $\*A$ of $A$ in our 
$\l^+$-saturated nonstandard universe $\*\V$.  

\begin{Th}   \label{nonstandard_appr}
The algebra $A$ is approximable by finite systems of a class $\K$ iff there exists a 
hyperfinitie system $H\in\*\K$ and an internal injection $j:H\to A$ that satisfy the following
properties.

\begin{enumerate}
\item $\forall a\in A\exists h\in H (j(h)\approx a)$.
\item $\forall h_1,h_2\in j^{-1}(\ns(\*A))\left(j(h_1\odot h_2)\approx j(h_1)\circ j(h_2)\right)$
\end{enumerate}
\end{Th}
 
{\bf Proof} $\Rightarrow$ Assume that $A$ is approximable by finite systems of a class $\K$. 
Let $\H$ be the set of all pairs $(C,\U)$, where $C\subseteq A$ is a compact and $\U$ --- a finite cover of $C$ by open sets. Notice that the cardinality of $\H$ is $\l$.
Consider the following preordering $\leq$ of $\H$:
$$
(C,\U)\leq(C',\U')\Leftrightarrow (C\supset C')\land\forall U\in\U\left(U\cap C'\neq\emptyset\Rightarrow\exists V\in\U'(U\subseteq V)\right) \eqno(3)
$$
Let us show that 
$$
\forall (C_1\U_1),\,(C_2,\U_2)\exists (C,\U)\left( (C,\U)\leq(C_1,\U_1)\land (C,\U)\leq(C_2,\U_2)\right).
\eqno (4)
$$
Let $C=C_1\cup C_2$ and $\U=\{U\setminus C_1\ |\ U\in\U_2\}\cup \{U\setminus C_2\ |\ U\in\U_1\}
\cup\{U\cap V\ |\ U\in\U_1,V\in\U_2,U\cap V\neq\emptyset\}.$

It is enough to show that $(C,\U)\leq (C_1,\U_1)$. By construction $C\supset C_1$ and if $W\in\U$
and $W\cap C_1\neq\emptyset$ then $W$ is either of the form $U\setminus C_2,\ U\in\U_1$ or of the form $U\cap V,\ U\in\U_1, V\in\U_2$. In both cases (3) holds.

Let $\H(C,\U)=\{(C',\U')\ |\ (C',\U')\leq (C,\U)\}$. The family 
$\{\H(C,\U)\ |\ (C,\U)\in\H\}$ is of the cardinality $\l$ and has 
the finite intersection property by (4). By $\l^+$-saturation
$$
\bigcap\limits_{(C,\U)\in\H}\*\H(C,\U)\neq\emptyset
$$
and thus there exists a pair $(C_0,\U_0)\in\*\H$ such that 
$\forall(C,\U)\in\H\ (C_0,\U_0)\leq
(\*C,\*\U)$. By definition of $\leq$ the compact $C_0\supset\bigcup\{\*C\ |\ C\subseteq A,\ C\ \mbox{is compact}\}=\ns(\*A)$ by Proposition~\ref{topprop}(4).

Let us show that 
$$
\forall a\in A\forall U\in\U_0(a\in U\Rightarrow U\subseteq\mu(a)) \eqno(5)
$$.

Indeed, let $V$ be any standard neighborhood of $a$. Consider $C=\{a\},\ \U=\{V\}$. Since
$(C_0,\U_0)\leq(\*C,\*\U)$, if $a\in U\in\U_0$ then by (3) $U\subseteq\*V$. Thus
$U\subseteq\bigcap\{\*V\ |V\ni a,\V\ \mbox{is open}\}=\mu(a)$.  

By the Transfer principle there exist a hyperfinite algebra $H$ and an internal injection
$j:H\to\*A$ such that $<H,j>$ is a $(C_0,\U_0)$-approximation of $\*A$. We have to show that 
$<H,j>$ satisfies the properties 1) and 2) of the theorem.

Let $a\in\ns(\*A)$.
By Definition~\ref{approx} and the Transfer principle since $a\in C_0$, 
$\exists U\in\U_0\exists
h\in H(a,j(h)\in U)$. By (5) $a\approx j(h)$ and thus the property 1) holds.
 
Let $j(h_1),j(h_2)\in\ns(\*A)$. 
Since $j(h_1),j(h_2)\in C_0$ and $j(h_1)\circ j(h_2)\in C_0$ 
(due to continuity of $\circ$) 
we have
$$
\exists U\in\U_0\ \left(j(h_1\odot h_2), j(h_1)\circ j(h_2)\in U\right).
$$
 Let us show that
$U\subseteq\mu(a)$ for some $a\in A$. This will prove the property 2).
Put $a=\o\left(j(h_1)\circ j(h_2)\right)$.  
It is enough to show that for an arbitrary open $V\ni a$ holds 
$U\subseteq\*V$. Consider an open relatively compact $W\ni a$ such that 
$\overline W\subseteq V$, where $\overline W$ is the closure of $W$. 
Such $W$ exists because a locally compact space is regular. 
Notice that $j(h_1)\circ j(h_2)\in \*W$ since 
$j(h_1)\circ j(h_2)\approx a$ by Proposition~\ref{topprop}. 
Let $C=\overline W$ and $\V=\{V\}$. We have $(C_0,\U_0)\leq (\*C,\*\V)$, 
$U\ni j(h_1)\circ j(h_2)$ and thus $U\cap\*C\neq\emptyset$. Now by (3) and 
the Transfer principle $U\subseteq\*V$.

$\Leftarrow$ Let $<H,j>$ satisfy the properties 1) and 2) of the theorem, 
$C\subseteq A$ 
be
a compact and $\U$ --- a finite cover of $C$ by open sets. 
We have to show that $j(H)$ is a $(\*C,\*\U)$-grid and
$j:H\to\*A$ is a $(\*C,\*\U)$-homomorphism. Then the theorem will be proved 
by Transfer principle
(working 
now the opposite direction of
the first part of the proof).

Let $\*C\cap\*U\neq\emptyset$ then, by Transfer principle 
$\exists c\in C\cap U$.
By the property 1) $\exists h\in H\,(j(h)\approx c)$. Thus $j(h)\in\*U$
by Proposition~\ref{topprop}.
This proves that $j(H)$ is a $(\*C,\*\U)$-grid.

Let $j(h_1), j(h_2), j(h_1)\circ j(h_2)\in \*C$ then by Proposition~\ref{topprop}(3)
$a_1=\o j(h_1),\ a_2=\o j(h_2),\ a= \o j(h_1)\circ j(h_2)\in C.$ By the property 2)
$j(h_1\odot h_2)\approx j(h_1)\circ j(h_2)\approx a$. 
Thus if $a\in U\in\U$ then $j(h_1\odot h_2),j(h_1)\circ j(h_2)\in\*U$ and, 
consequently, $<H,j>$ is a 
$(\*C,\*\U)$-homomorphism. $\Box$

\begin{Prop} \label{nonst-gr-appr}
A locally compact group $G$ is approximable by finite algebras of a class $\K$ in the sense of Definition~\ref{gr_approx} iff hyperfinitie system $H\in\*\K$ and an internal injection $j:H\to A$ that satisfy the properties 1) and 2) of Theorem~\ref{nonstandard_appr}.
\end{Prop}

Proposition~\ref{equiv-for-groups} follows immediately from Theorem~\ref{nonstandard_appr}
and Proposition~\ref{nonst-gr-appr}.

{\bf Proof}. The proof of Proposition~\ref{nonst-gr-appr} is very similar to the proof of 
Theorem~\ref{nonstandard_appr} but simpler, so we'll sketch briefly the main points.

$\Rightarrow$ By the $\l$-saturation of the nonstandard universe there exist a compact set $C\subseteq\*G$ such that $\ns(\*G)\subseteq C$ and an open neighborhood 
$U\subseteq\*G$ of the unit $e$ such that 
$U\subseteq\mu(e)$. By Transfer principle there exists a hyperfinite $(C,U)$-approximation 
$<H,j>$ of $\*G$ such that $H\in\*K$. It is easy to see that $<H,j>$ satisfies the properties 1) and 2) of Theorem~\ref{nonstandard_appr}.

$\Leftarrow$ Let $<H,j>$ be a hyperfinite algebra that satisfies the properties 1) and 2) of Theorem~\ref{nonstandard_appr} and $H\in\*K$. Then it is easy to see that for any standard $C\subseteq G$ and any neighborhood $U\subseteq G$ of the unit $<H,j>$ is a $(\*C,\*U)$-approximation of $\*K$. 
Now by Transfer principle the condition 4) of Definition~\ref{gr_approx} $\Box$
\section{Proof of Theorem~\ref{functional}}

We start with construction of a functional $I$ that satisfies the condition of Theorem~\ref{functional}.
We'll consider only the case of a locally compact $l$-quasigroup $A$, approximable by finite $l$-quasigroups. The case of $r$-quasigroups is similar.

Let $<H,j>$ be a hyperfinite $l$-quasigroup that satisfies Theorem~\ref{nonstandard_appr}.
Let $V\subseteq A$ be a compact set with the nonempty interior i.e.
there exists a nonempty open set $U\subseteq V$ and thus, by the regularity of the topological space
$A$ there exists an open $W$ such that $\overline W\subseteq U$. 

We'll write $W\sqsubset D$ if $\overline W\subseteq U\subseteq D$ for an open set
$U$.
In the proofs we often will use the following\\
{\bf Statement.} If $D$ is a compact, $W\sqsubset D$, 
$x\in \* W$ and $y\approx x$ then $y\in\* D$. \\
The statement easily follows from Proposition~\ref{topprop}.

Let $\Delta^{-1}=|j^{-1}(\* V)|$.
Define the functional $I(f)$ as follows:
$$
I(f)=\o \left(\Delta \sum_{h\in H}\* f(j(h))\right) \eqno(6)
$$

The proof of Theorem~\ref{functional} follows from two following lemmas: Lemma~\ref{Lm.cont} and Lemma~\ref{Lm.dec}.

\begin{Lm} \label{Lm.cont}
The functional $I(\cdot )$ is a Radon measure on
$C_0(A)$.
\end{Lm}

We need three following propositions.

\begin{Prop}\label{Prop.1}
Let $D\subseteq A$ be compact and $U\sqsubset D$ be an open set.
Then for all $a\in A$ the following inequality holds
$$
|j^{-1}(a\circ\*U)|\leq |j^{-1}(\*D)|
$$
\end{Prop}
{\bf Proof}. 
Let $x\in j^{-1}(a\circ\*U)$, or $j(x)\in a\circ\*U\subseteq\ns(\*A)$. By the left cancellation law and   Transfer principle
$b_a\circ j(x)\in\*U$, where $b_a\in A$  does not depend on $x$. 

By Theorem~\ref{nonstandard_appr} there exists $\beta\in H$, such that
$b_a\approx j(\beta)$. 
So, $b_a\circ j(x)\approx j(\beta\circ x)\in \*D$ because $U\sqsubset D$. 
Consequently, 
$\beta\circ (j^{-1}(a\circ\*U))\subseteq j^{-1}(\*D)$, but the function 
$l_\beta(x)=\beta\circ x$ is an
injection, because $H$ is an l-quasigroup. $\Box$.

\begin{Prop}\label{Prop.2}
Let $X, Y\subseteq A$ be compact sets and $Y$ has the nonempty interior. 
Then there exists $0<C_{X,Y}\in \R$, such that
$$
{|j^{-1}(\*X)|\over |j^{-1}(\*Y)|}\leq C_{X,Y}
$$
\end{Prop}
{\bf Proof.} Take an open $U\sqsubset Y$. By the definition of a 
topological $l$-quasigroup  
the mapping
$l_a:A\to A$ is a continuous homeomorphism for any $a\in A$.  
Thus $l_a(U)=a\circ U$ 
is an open set for any $a\in A$. By the definition of an $l$-quasigroup $\forall z\in A\ A\circ z=A$. Thus $A\circ U$ covers $X$ and so, there exists a finite set $F\subseteq A$ such that
$X\subseteq F\circ U$. It means that $\*X\subseteq F\circ\*U$ ($\*F=F$).
Consequently,
$$
|j^{-1}(\*X)|\leq \sum_{\alpha\in F}|j^{-1}(\alpha\circ \*U)|,
$$
and, by Proposition~\ref{Prop.1} $|j^{-1}(\*X)|\leq |F|\cdot |j^{-1}(\*Y)|$.
So, one can take $C_{X,Y}=|F|$. $\Box$.
\begin{Prop} \label{Prop.3}
Let $\phi:H\to \*\R$ satisfy the following conditions: 
\begin{enumerate}
\item $\forall h\in H\;\;\f(h)\geq 0$;
\item $j(\supp (\f))\subseteq \*S$, where $S\subseteq A$ is a compact;
\item If $j(h)\in\* D$, then $\f(h)>\alpha$ for some compact $D\subseteq A$ with the
nonempty interior and some $\alpha\in \*\R$.
\end{enumerate}
Then
$$
{1\over C_{V,D}}\alpha\leq \Delta\sum_{h\in H} \phi(h)\leq C_{S,V}
\sup(\phi). \eqno(7)
$$
\end{Prop}
The proof of Lemma~\ref{Lm.cont}  follows immediately from 
Proposition~\ref{Prop.3}.

Indeed, take $\f(h)=\*f(j(h))$ for any $0< f\in\C_0(A)$. Then $\f$ satisfies the conditions of 
Proposition~\ref{Prop.3}. Obviously $\f(h)\geq 0$, $S=\supp(f)$. Since $f>0$ there exists a point
$a\in A$ such that $f(a)>0$ and thus there exist an open $U\ni a$ and a positive $\a$ such that
$\forall b\in U\ f(b)>\a$. Take any relatively compact open $W$ such that $\overline W\subseteq U$.
Then $D=\overline W$ satisfies the condition 3 of Proposition~\ref{Prop.3}. 

By (6) and the first inequality (7) $I(f)\neq 0$. By the second inequality (7) the linear functional $I$
is continuous. 

{\bf Proof of Proposition~\ref{Prop.3}}. Recall that 
$\Delta^{-1}=|j^{-1}(\* V)|$.
By Proposition~\ref{Prop.2}:
$$
 \Delta\sum_{h\in H} \phi(h)\geq  \Delta\sum_{j(h)\in \*D} \phi(h)
\geq {\alpha\over C_{V,D}}.
$$
This proves the first inequality (7).
The second inequality (7) is obtained as follows:
$$
 \Delta\sum_{h\in H} \phi(h) = \Delta\sum_{j(h)\in \* S} \phi(h)
\leq C_{S,V}\cdot \sup_x\phi(x)\ \ \ \Box
$$ 
\begin{Lm} \label{Lm.dec}
The functional $I$, defined by (6), satisfies the inequality 
$I(f)\geq I(l_a(f))$
for any non-negative $f\in C_0(A)$.
\end{Lm}
{\bf Proof}.
For $X\subseteq A$, $z\in A$ let us denote 
$/(X,z)=\{/(x,z)\; :\; x\in X\}$, see Definition~\ref{quasigroups}.
Let $S\subseteq A$ be a compact, 
such that there exists 
an
open set $U\subseteq A$ with the property
$/(\supp (f),a)\subseteq U\sqsubset S$. 
Let $h\in H$ such that $j(h)\approx a\in A$.
First of all the following equality holds:
$$
\o\left(\Delta\sum_{x\in H}\* f(a\circ j(x))-
\Delta\sum_{x\in H}\* f(j(h)\circ j(x))\right)=0   \eqno(8)
$$
To prove it let
$\f(x)=|\* f(a\circ j(x))-\* f(j(h)\circ j(x))|$ and apply
Proposition~\ref{Prop.3} as follows. 
By the continuity of $\circ$ and $/(\cdot,\cdot)$ in $A$, we have

$a\circ j(x)\in\ns\Leftrightarrow j(x)\in\ns\Leftrightarrow j(h)\circ j(x)\in\ns$. 
Let us show now that 
$j(\supp(\f))\subseteq\*S$. 

It is enough to show that 
$$
j(x)\notin\* S\Longrightarrow \*f(a\circ j(x))=\*f(j(h)\circ j(x))=0 \eqno (9)
$$ 

Assume, that $a\circ j(x)\in \*\supp(f)$ and thus 
$j(x)\in/(\*\supp(f),a)\subseteq \*S$. This proves the first equality. 

Assume, that $j(h)\circ j(x)\in\*\supp(f)$. 
Then $j(x)\in/(\*\supp(f),j(h))\approx/(\*\supp(f),a)\subseteq \*U$. 
But $U\sqsubset S$ and $j(x)\in \*S$.
So, we get the contradiction.

Since $a\circ j(x)\approx j(h)\circ j(x)$ if $j(x)\in\ns(\* A)$ and $\supp (\f)\in j^{-1}(\* S)\subseteq j^{-1}
(\ns(\*A))$ we 
have
$\supp(\f)\approx 0$ and by the second inequality (7) 
$\Delta\sum\limits_{h\in H}\f(h)\approx 0$. This proves (8).

Let us show that the following inequality holds
$$
\o\left(\Delta\sum_{x\in H}\* f(j(h\circ x))-
\Delta\sum_{x\in H}\* f(j(h)\circ j(x))\right)\geq 0   \eqno(10)
$$
By (9) we have
$$
\Delta\sum_{x\in H}\* f(j(h\circ x))-
\Delta\sum_{x\in H}\* f(j(h)\circ j(x))=
\Delta\sum_{j(x)\not\in\* S}\* f(j(h\circ x))+
\Delta\sum_{j(x)\in\* S}(\* f(j(h\circ x))- \* f(j(h)\circ j(x)))
$$
Obviously,
$$
\Delta\sum_{j(x)\not\in\*S}\* f(j(h\circ x))=c\geq 0.
$$
But
$$
\Delta\sum_{j(x)\in\* S}(\* f(j(h\circ x))- \* f(j(h)\circ j(x)))\approx 0.
$$
Indeed, since $j(h),j(x)\in\ns(\* A)$, when $j(x)\in\* S$ we have 
$j(h\circ x)\approx j(h)\circ j(h)$
by Theorem~\ref{nonstandard_appr} and thus $\* f(j(h\circ x))\approx \* f(j(h)\circ j(h))$ by the continuity of $f$. Thus $\b=\sup\limits_{j(x)\in\* S}|\* f(j(h\circ x))- \* f(j(h)\circ j(x))|)\approx 0$
and by Proposition~\ref{Prop.2}
$$
|\Delta\sum_{j(x)\in\* S}(\* f(j(h\circ x))- \* f(j(h)\circ j(x)))|\leq C_{S,V}\b\approx 0
$$

Since $\{h\circ x\ |\ x\in H\}$ is a permutation of $H$ we have 
$$
\Delta \sum_{x\in H}\* f(j(x))=\Delta \sum_{x\in H}\* f(j(h\circ x))
$$
Now
$$
I(f)-I(l_a(f))=\o\left(\Delta \sum_{x\in H}\* f(j(x))-\Delta \sum_{x\in H}\* f(a\circ j(x))\right)=
$$
$$
=\o\left(\left(\Delta \sum_{x\in H}\* f(j(h\circ x))-\Delta \sum_{x\in H}\* f(j(h)\circ j(x))\right)+\left(\Delta \sum_{x\in H}\* f(j(h)\circ j(x))-\Delta \sum_{x\in H}\* f(a\circ j(x))\right)\right)
$$

The first term on the RHS of this equality is positive by (10), the second --- infinitesimal by (8)
and so $I(f)-I(l_a(f))\geq 0$ $\Box$

\section{Proof of Theorem~\ref{semigroups}}

First of all we'll formulate some necessary results about the structure of finite semigroups from \cite{RT}, where one can find the proofs.

Let $S$ be a finite semigroup.
\begin{Def} \label{fin.sem.}
\begin{enumerate}
\item
$x\in S$ is called to be zero ($x=0$) iff $\forall y\in S \;xy=yx=x$.
(Obviously if the zero exists it is unique).
\item
The set $I\subseteq S$ is a left (right) ideal iff 
$SI\subseteq I$ ($IS\subseteq I$). $I\subseteq S$ is an ideal iff $I$ 
is a left and a right ideal.
(Obviously an ideal (a left or a right ideal) is a subsemigroup.)
\item
$S$ is called to be $0$-simple if it has no proper ideals but $\{ 0\}$ and
$\emptyset$.
\item
$S$ is a zero semigroup iff $\forall s,t\in S \;st=0$
\item
Let $I\subseteq S$ be an ideal of $(S,\cdot)$. The quotient semigroup
$S/I$ is the set $(S\backslash I)\cup\{ 0 \}$ with multiplication 
"$*$" defined as the follows
$$
s_1* s_2=\left\{                \begin{array}{l}			
		s_1\cdot s_2,\; \mbox{if}\; s_1\cdot s_2\notin I\\
		0,\; \mbox{if}\; s_1\cdot s_2 \in I	
	                         \end{array}		   
               \right.
$$
\item
Maximal sequence of ideals for $S$ is the ordered sequence of ideals of $S$
$$
S=I_0\supset I_1\supset I_2... I_n\supset I_{n+1}=\emptyset,
$$
such that there are no ideals $I'$ of $S$, $I_k\supset I'\supset I_{k+1}$.
\end{enumerate}
\end{Def}

It is clear that any finite semigroup has a maximal sequence of ideals.

\begin{Th} \label{T_simp}
Any semigroup $I_{r-1}/I_r$ is  $0$-simple or zero.
\end{Th}
Let $n,m\in \N$, $H$ be a group,
$\rho:\{1,..,n\}\times\{1,..,m\}\to H\cup\{ 0\}$. 
Consider Rees semigroup
$$
S(n,m,H,\rho)=\{(i,j,h),\;\; i=1,...,n;\;j=1,...,m;\;h\in H\} \cup \{ 0\}
$$
$$
(i_1,j_1,h_1)(i_2,j_2,h_2)=
	\left\{	\begin{array}{l}			 (
           i_1,j_2, h_1\rho(i_2,j_1)h_2),\;\mbox{if}\;\rho(i_2,j_1)\in H\\
			 0,\;\mbox{if}\; \rho(i_2,j_1)=0	     
                  \end{array}	   
          \right.
$$
Rees semigroup is called regular if $\forall i\exists j\;\rho(i,j)\neq 0$ 
and $\forall j\exists i\;\rho(i,j)\neq 0$.

\begin{Th} \label{structure}
Any finite $0$-simple semigroup $S$ (with zero) is isomorphic to a regular 
Rees semigroup.
\end{Th}
(If $S$ is a semigroup without zero we may add zero to $S$ or remove zero
from the Rees semigroup.) \\
This theorem implies
\begin{Cor}  \label{C_left}
Let $S$ be a $0$-simple finite semigroup, $0\neq s\in S$ and $F=sSs$.
Then $F$ is a zero subsemigroup or, $F\backslash \{0\}$ is a group.
\end{Cor}
{\bf Proof} Let $s=(i_s,j_s,h_s)$. If $F\neq \{0\}$, then
$F=\{sas,\; a\in S\}=\{(i_s, j_s, h),\;h\in H\} \cup \{0\}$.
If $\rho(i_s,j_s)=0$, then $F$ is a zero semigroup; if $\rho(i_s, j_s)=g$,
then $\phi :F\backslash \{0\}\to H$, $\phi(i_s,j_s,h)=hg$, is an isomorphism.

We are able now to prove Theorem~\ref{semigroups}.

Let $G$ be a locally compact group that is approximable by finite semigroups and $<S,\phi>$ --- a hyperfinite approximation of $G$ by a hyperfinite semigroup $S$ that exists by Theorem~\ref{nonstandard_appr}. We denote both
operations --- in $G$ and in $S$ by $\cdot$ since this does not lead to any misunderstanding.

Consider an internal hyperfinite maximal sequence of ideals in $S$
(see definition~\ref{fin.sem.}(6)) $$
S=I_0\supset I_1\supset I_2 ...I_n\supset I_{n+1}=\emptyset .
$$
that exists by Transfer principle.

By assumptions, $\phi(S)\cap \ns\neq\emptyset$ and, consequently,
there exists $r\in\*\N$ such that $\phi(I_{r-1})\cap \ns(\* G)\neq\emptyset$ and
$\phi(I_r)\cap\ns(\*G)=\emptyset$. There are two cases.
\begin{enumerate}
\item $I_r=\emptyset$. 
Then take $F=I_{r-1}$ and $\psi=\phi|_F$.

\item $I_r\neq\emptyset$. 
Then take $F=I_{r-1}/I_r$ and $\psi:F\to G$, $\psi(x)=\phi(x)$, for 
$x\neq 0$ and $\psi(0)=g\notin\ns(\*G)\cup\mbox{Im}(\phi)$. 
Such $g$ exists by the following reasons. The group $G$ is not compact, 
otherwise $\*G=\ns(\*G)$, but
$\phi(I_r)\cap\ns(\*G)=\emptyset$, $I_r\neq\emptyset$
 and $\phi(S)\subseteq G$. It is easy to see that there exist an internal compact $D\supset\ns(\*G)$. The set
$\*G\setminus\*D$ is not compact and thus not hyperfinite. So $\*G\setminus\left(\ns(\*G)\cup \mbox{Im}(\phi)\right)\neq \emptyset$.
\end{enumerate}
Let us 
prove that $<F,\psi>$ approximates $G$ in the sense of 
Theorem~\ref{nonstandard_appr}. 
Let us denote by $*$ the operation on $F$.

First, we will show that $\psi$ is an almost homomorphism.
Let $x, y \in F$ and $\psi(x),\psi(y)\in\ns(\*G)$, we have to prove that 
$\psi(x*y)\approx \psi(x)\psi(y)$. For the case 1) it is trivial, since $\psi$
is a restriction of $\phi$ on subsemigroup. Consider case 2). 
Since $\psi(x),\psi(y)\in\ns(\*G)$, one has $x, y\neq 0$, so, 
$\psi(x)=\phi(x)$ and $\psi(y)=\phi(y)$. Then 
$\ns(\*G)\ni\phi(x)\phi(y)\approx\phi(xy)$. So, $\phi(xy)\in\ns$ and thus 
$xy\not\in I_r$. By the definition  of the operation in a quotient semigroup   
$x*y=xy\neq 0$, and by the construction $\psi(x*y)=\phi(xy)$.
 
It remains to prove that
$\forall g\in G\exists x\in F\; g\approx\psi(x)$ or, the same,
$\forall g\in G\exists x\in I_{r-1}\; g\approx\phi(x)$.
Since $\phi(I_{r-1})\cap\ns(\*G)\neq\emptyset$, there exists an
$x\in I_{r-1}$ such that $\phi(x)\in\ns(\*G)$. Since $\mbox{}^{-1}$ is a continuous operation
$(\phi(x))^{-1}\in\ns(\*G)$ and there exists $y\in S$ 
$\phi(y)\approx(\phi(x))^{-1}$. So, $e\approx\phi(y)\phi(x)\approx\phi(yx)$. Notice that $yx\in I_{r-1}$ 
since $I_{r-1}$ is an ideal and $x\in I_{r-1}$. Now, let $g\in G$ and $s\in S$ such that 
$\phi(s)\approx g$. Then $\phi(yxs)\approx g$ and $yxs\in I_{r-1}$.

Obviously a zero semigroup can never approximate
an infinite group, and thus $F$ is $0$-simple hyperfinite semigroup by Theorem~\ref{T_simp}

Let $s\in F$ be such that $\psi(s)\approx e$. Consider the semigroup $T=sFs$ It is easy to see if $j=\psi |_T$ then the pair $<T,j>$ approximates $G$.

By Corollary~\ref{C_left} of Theorem~\ref{structure}
$H=T\setminus\{0\}$, is a hyperfinite group. If $T$ does not contain $0$,
then the proof is done. Suppose $0\in T$. Then it is enough to prove that $j(0)\not\in\ns(\*G)$. 

Suppose that $j(0)\in\ns(\*G)$.
If $j(0)\approx e$, then $\forall x\in G\; xe\approx e$ which is impossible.
If $j(0)\approx x$ and $x\neq e$, then there exist $y$ such that $j(y)\approx x^{-1}$. Now\\
$e=xx^{-1}\approx j(0)j(y)\approx j(0*y)=j(0)\approx x$.
This is impossible since $x,e$ are standard and $x\neq e$. $\Box$

\bigskip

Instituto de Investigacion en Communicacion Optica de Universidad Autonoma de 
San Luis Potosi, Mexico

Eastern Illinois University, USA

\bigskip

IICO-UASLP

AvKarakorum 1470

Lomas 4ta Session

San Luis Potosi SLP 7820

Mexico

Phone: 52-444-825-0892 (ext. 120)

e-mail:glebsky@cactus.iico.uaslp.mx\\

\bigskip

Mathematics Department 1-00036

Eastern Illinois University

600 Lincoln Avenue

Charleston, IL 61920-3099

USA

Phone: 1-217-581-6282

e-mail: cfyig@eiu.edu 

\bigskip

1991 {\it Mathematics Subject Classification}. Primary 26E35, 03H05; Secondary 28E05, 42A38


\begin{thebibliography}{100}
\bibitem{Gor} Gordon E. Nonstandard Methods in Commutative Harmonic Analysis.
AMS, Providence, Rhode Island, 1997.
\bibitem{VG} A.M.Vershik, E.I.Gordon. Groups locally embedded into the class of finite groups. (Russian) Algebra i Analiz 9 (1997),no. 1, p. 71-97; translation in St.Petersburg Math. J. 9 (1998), no. 1, 49-67.
\bibitem{AGKh} Albeverio S., Gordon E., Khrennikov A. 
Finite dimensional approximations of operators in the spaces of functions on
locally compact abelian groups. {\it Acta Applicandae Mathematicae} {\bf} 64(1) pp. 33-73, October 2000
\bibitem{AGG} M.A.Alekseev, L.Yu. Glebskii, E.I.Gordon. On approximations of groups, group actions and Hopf algebras. {\it Representation Theory, Dynamical Systems,
Combinatorial ana Algebraic Methods. III}, A.M.Vershik – editor,
Russian Academy of Science. St.Petersburg Branch of V.A.Steklov's Mathematical
Institute. Zapiski nauchnih seminarov POMI  256 (1999), 224-262. (in Russian; Engl. Transl. in {\it Journal of Mathematical Sciences}, 107, No.5 (2001), pp.4305-4332)
\bibitem{Gr} M.Gromov Endomorphisms of symbolic algebraic varieties. {\it J.Eur.Math. Soc}. 1 (1999), 109 – 197
\bibitem{Resv} E.I. Gordon, O.A. Rezvova. On hyperfinite approximations of
the field $\R$.{\it Rueniting the Antipodes --- Constructive and Nonstandard Views of the Continuum,
Proceedings of the Symposium in San Servolo/Venice, Italy, May 17 –20, 2000}.
B.Ulrich, H.Ossvald and P. Schuster, editors. Synth\'eseLibrary, volume 306
By Kluwer Academic Publishers, Dordrecht, etc., 2001
\bibitem{Alb} Albeverio S., Fenstad J-E., Hoegh-Krohn R.,
Lindstroem T. {\it Nonstandard Methods in Stochastic Analysis and
Mathematical Physics}. Academic Press, New York,1986
\bibitem{Loeb} Nonstandard Analysis for the Working Mathematicians. P.A. Loeb and M.P.H. Wolff,
editors. Mathematics and Applications, volume 510. Kluwer Academic Publishers, Dordrecht/Boston/London, 2000 
\bibitem{quas} Quasigroups and Loops. Theory and Applications. O. Chein, H.O. Pfulgfelder and
J.D.H. Smith, editors. Sigma Series in Pure Mathematica, volume 8. Heldermann Verlag, Berlin, 1990.
\bibitem{Ryser} H. J. Ryser, Combinatorial mathematics, The Carus
Mathematical Monographs, 15, The Mathematical Association of America, 1963.
\bibitem{Halmos} Paul R. Halmos, Measure Theory, Springer-Verlag, New York
1974.
\bibitem{Ross} D.A. Ross Measures invariant under local homeomorphisms {\it Proc. Amer. Math. Soc.}
{\bf 102}(4), 1988, pp 901-905 
\bibitem{Bos} B.Zivaljevic. Hyperfinite transversal theory. {\it Ph.D. Thesis. Univ. of Illinois at
Urbana - Champaign}, 1989
\bibitem{RT} J. Rhodes, B. Tilson, Theorems on local structure of 
finite semigroup, {\it in Algebraic theory of machines, languages and 
semigroups,ed. M.A. Arbib} Acad. Press, New York \& London, 1968.
\end{thebibliography}
\end{document}